\documentclass[english,conference]{ieeeconf}
\usepackage[T1]{fontenc}
\usepackage[latin9]{inputenc}
\usepackage{mathtools}
\usepackage{amsmath}
\usepackage{amsthm}
\usepackage{amssymb}
\usepackage{graphicx}

\makeatletter
\theoremstyle{plain}
\newtheorem{thm}{\protect\theoremname}
\theoremstyle{plain}
\newtheorem{prop}[thm]{\protect\propositionname}
\theoremstyle{plain}
\newtheorem{lem}[thm]{\protect\lemmaname}

\usepackage[caption=false,font=footnotesize]{subfig}
\usepackage{url}
\usepackage[nocompress]{cite}

\usepackage{amsmath}
\usepackage{amsfonts}

\newcommand{\norm}[1]{\left\Vert#1\right\Vert}

\newcommand{\optmin}{\mathrm{min.}}

\newcommand{\matJ}{J}
\newcommand{\matdJ}{J_\perp}

\newcommand{\vOne}{\mathbf{1}}

\newcommand{\opt}{^{\star}}

\newcommand{\xmin}{x_0\opt}
\newcommand{\xopt}{x^{\star}}
\newcommand{\uopt}{u^{\star}}
\newcommand{\xiopt}{\xi^{\star}}
\newcommand{\fopt}{f_0^{\star}}

\newcommand{\deq}{\coloneqq}
\newcommand{\bigO}{O}

\newcommand{\sr}{\sigma}

\usepackage{amsthm}
\theoremstyle{plain}
\@ifundefined{thm}{%
\newtheorem{thm}{Theorem}%
}{}
\newtheorem{assumption}[thm]{Assumption}

\allowdisplaybreaks[3]
\IEEEoverridecommandlockouts
\newcommand{\allthanks}{\thanks{The author is with the Department of Electrical and Computer Engineering, University of Illinois, Chicago, IL 60607. {\tt\scriptsize {hanshuo@uic.edu}}.}}

\makeatother

\usepackage{babel}
\providecommand{\lemmaname}{Lemma}
\providecommand{\propositionname}{Proposition}
\providecommand{\theoremname}{Theorem}

\begin{document}

\title{\textbf{\Large{}Computational Convergence Analysis of Distributed Gradient Tracking\\ for Smooth Convex Optimization Using Dissipativity Theory}}

\author{Shuo Han\allthanks}
\maketitle
\begin{abstract}
We present a computational analysis that establishes the $\bigO(1/K)$ convergence of the distributed gradient tracking method when the objective function is smooth and convex but not strongly convex. The analysis is inspired by recent work on applying dissipativity theory to the analysis of centralized optimization algorithms, in which convergence is proved by searching for a numerical certificate consisting of a storage function and a supply rate. We derive a base supply rate that can be used to analyze distributed optimization with non-strongly convex objective functions. The base supply rate is then used to create a class of supply rates by combining with integral quadratic constraints. Provided that the class of supply rates is rich enough, a numerical certificate of convergence can be automatically generated following a standard procedure that involves solving a linear matrix inequality. Our computational analysis is found capable of certifying convergence under a broader range of step sizes than what is given by the original analytic result. 

\end{abstract}

\section{Introduction}

Distributed optimization algorithms have a wide range of applications in engineering~\cite{molzahn_survey_2017,nedic_distributed_2018,rabbat_distributed_2004} and statistics~\cite{boyd_distributed_2011} when the scale of the optimization problem becomes too large to be solved centrally. A fundamental issue in the analysis of optimization algorithms is convergence, in particular convergence rate, which is a measure of how quickly an algorithm is able to locate an optimal solution. Traditional analysis of convergence rates relies on nonconstructive analytic proof techniques, which are often devised on an algorithm-by-algorithm basis and therefore do not readily generalize to new algorithms. As a result, one often needs to start the analysis from scratch when new requirements such as robustness, security, and communication constraints are introduced to existing algorithms. 

Recently, there has been work on using computational tools for analyzing the convergence rate of optimization algorithms. This computational approach is analogous to searching for a Lyapunov function for nonlinear systems using the sums-of-squares technique~\cite{parrilo_structured_2000}, in which a numerical certificate is used to prove desired properties. Our work in this paper is heavily inspired by the work of Lessard et al.~\cite{lessard_analysis_2016}, in which an optimization algorithm is viewed as a feedback interconnection of a linear dynamical system and a nonlinear memoryless uncertain system. The linear system captures the update rules, whereas the nonlinear system captures properties of the objective function such as strong convexity and smoothness. The nonlinear system is treated as uncertain because  convergence of optimization algorithms often needs to be established over a class of objective functions as opposed to a specific one. From this viewpoint, convergence of optimization algorithms becomes equivalent to stability of the feedback interconnection, which can be analyzed using tools from nonlinear and robust control. In many cases, the nonlinear uncertain system may be described by \emph{integral quadratic constraints} (IQCs), first proposed by Megretski and Rantzer~\cite{megretski_system_1997}. Because IQCs only rely on properties of the objective function, they only need to be derived once and can be reused in the analysis of new optimization algorithms. For strongly convex objective functions, the method in~\cite{lessard_analysis_2016} has been shown to prove exponential (linear in the language of optimization theory) convergence for a variety of algorithms such the gradient descent method and Nesterov's method. The computational analysis is automated and only requires solving a linear matrix inequality (LMI), which can be computed efficiently using existing optimization software. 

There have been two important extensions of~\cite{lessard_analysis_2016}. One extension handles non-strongly convex objective functions, which arise from a number of applications in machine learning~\cite{bach_non-strongly-convex_2013}. Without the assumption of strong convexity, optimization algorithms typically exhibit sublinear convergence, which cannot be certified by the quadratic Lyapunov-like functions used in~\cite{lessard_analysis_2016}. Fazlyab et al.~\cite{fazlyab_analysis_2018} use a non-quadratic and time-varying Lyapunov function and, by combining with IQCs, prove sublinear convergence for a number of first-order methods (in a centralized setting). Alternatively, Hu and Lessard~\cite{hu_dissipativity_2017} use dissipativity theory to show sublinear convergence by searching for a suitable storage function. The other extension handles optimization algorithms in a distributed setting, in particular, distributed first-order methods for consensus optimization~\cite{sundararajan_robust_2017}. The challenge in distributed optimization is that the optimality condition is more complex than that in centralized optimization. As a result, some of the IQCs derived in the centralized setting cannot be applied directly for analyzing distributed algorithms. In the joint setting of distributed optimization and non-strongly convex objective functions, the best convergence rate has been shown to be $\bigO(1/K)$ by Qu and Li\cite{qu_harnessing_2018} and is achieved by an algorithm based on gradient tracking. However, the proof therein is based entirely on analytic techniques and relies on carefully bounding various quantities computed from the optimization algorithm. 

\paragraph*{Contribution}

In this paper, we give a computational analysis on the $\bigO(1/K)$ convergence of the distributed gradient tracking algorithm in~\cite{qu_harnessing_2018} when the objective function is smooth but non-strongly convex. Our main contribution is the derivation of a base supply rate for this setting, after which the standard procedure in~\cite{lessard_analysis_2016} can be used to compute a certificate of convergence. Compared to the analytic approach in~\cite{qu_harnessing_2018}, our computational approach is found to be able to certify convergence under a broader range of step sizes.

\section{Problem Description~\label{sec:problem}}

\subsection{Notation}

Denote by $\vOne$ the column vector of all ones, $I_{n}$ the $n\times n$ identity matrix (size omitted when clear from the context), and $\norm{\cdot}$ the $\ell_{2}$-norm of a vector. For a symmetric matrix $P$, we write $P\succeq0$ if $P$ is positive semidefinite. For a differentiable function $f$, we denote by $\nabla f$ the gradient of $f$. We reserve the subscript for indexing the entries of a vector and the superscript for indexing a given sequence (of either scalars or vectors). For example, the $i$-th entry of a vector $x\in\mathbb{R}^{n}$ is denoted by $x_{i}$, whereas a sequence of vectors is denoted by $\{x^{k}\}_{k\geq0}\coloneqq\{x^{0},x^{1},\dots\}$. 

\subsection{Problem description}

We consider an optimization problem of the form
\begin{equation}
\underset{x_{0}\in\mathbb{R}^{d}}{\optmin}\quad f_{0}(x_{0})\coloneqq\frac{1}{n}\sum_{i=1}^{n}f_{i}(x_{0}),\label{eq:prob}
\end{equation}
where $f_{i}\colon\mathbb{R}^{d}\to\mathbb{R}$ is convex for all $i\in\{1,2,\dots,n\}$. The set $\mathcal{X}_{\min}\coloneqq\mathrm{argmin}_{x_{0}\in\mathbb{R}^{d}}f_{0}(x_{0})$ of minimizers is assumed to be nonempty. As a result, the optimal value of the problem is well defined and is denoted by $\fopt\coloneqq\min_{x_{0}\in\mathbb{R}^{d}}f_{0}(x_{0})$. Throughout this paper, we make the following assumption on $f_{i}$.

\begin{assumption}\label{asu:smoothness}For $i\in\{1,2,\dots,n\}$, the function $f_{i}$ is \emph{$\beta$-smooth}, i.e., there exists $\beta>0$ such that
\[
\norm{\nabla f_{i}(x)-\nabla f_{i}(y)}\leq\beta\norm{x-y}
\]
holds for all $x,y\in\mathbb{R}^{d}$. 

\end{assumption}Note that we do not make the assumption that $f_{i}$ is strongly convex. To simplify notation, we derive our result only for $d=1$, and we will show in Section~\ref{subsec:Generalizations} that the main result (Theorem~\ref{thm:main}) also holds for $d>1$.

An algorithm is said to solve the optimization problem~(\ref{eq:prob}) if it can find at least one $x_{0}\opt\in\mathcal{X}_{\min}$. In this paper, we study the distributed gradient tracking algorithm in~\cite{qu_harnessing_2018} or equivalently the DIGing algorithm in~\cite{nedic_achieving_2017}, in which the update equations are given by\begin{subequations}\label{eq:dgd}
\begin{align}
x^{k+1} & =Wx^{k}-\eta s^{k}\label{eq:dgd-a}\\
s^{k+1} & =Ws^{k}+g(x^{k+1})-g(x^{k}),\label{eq:dgd-b}
\end{align}
\end{subequations}where 
\begin{gather}
x^{k}\coloneqq\left[\begin{array}{cccc}
x_{1}^{k} & x_{2}^{k} & \cdots & x_{n}^{k}\end{array}\right]^{T}\in\mathbb{R}^{n}\label{eq:xk}\\
g(x^{k})\coloneqq\left[\begin{array}{cccc}
\nabla f_{1}(x_{1}^{k}) & \nabla f_{2}(x_{2}^{k}) & \cdots & \nabla f_{n}(x_{n}^{k})\end{array}\right]^{T}\in\mathbb{R}^{n},\label{eq:gk}
\end{gather}
and $\eta>0$ is the step size. The initial condition is given by $s^{0}=g(x^{0})$, and $x^{0}$ is arbitrary. The matrix $W$ is assumed to be doubly stochastic (i.e., $W\vOne=\vOne$ and $\vOne^{T}W=\vOne^{T}$) and irreducible, which is a common assumption in the literature. It has been shown in~\cite[Th.~3]{qu_harnessing_2018} that the distributed gradient tracking algorithm is able to find a solution of problem~(\ref{eq:prob}) with an $\bigO(1/K)$ convergence rate; i.e., 
\[
\frac{1}{n}\sum_{i=1}^{n}\left[f_{0}(\hat{x}_{i}^{K})-\fopt\right]\leq\frac{V_{0}}{K+1},
\]
where $\hat{x}_{i}^{K}\deq\frac{1}{K+1}\sum_{k=0}^{K}x_{i}^{k}$, holds for some constant $V_{0}$ that depends on $W$, $\eta$, and the initial condition but not on $K$.

In this paper, we present a numerical procedure that proves the $\bigO(1/K)$ convergence rate using dissipativity theory. Compared with the result in~\cite{qu_harnessing_2018}, our procedure is able to certify convergence under a broader range of step sizes $\eta$. (See Section~\ref{sec:discussions} for a detailed comparison.) The conservatism of the result in~\cite{qu_harnessing_2018} is due to the fact that the proof of convergence therein relies on bounding various quantities computed from the optimization algorithm, where some bounds are made loose on purpose to keep the expressions in closed form. While some conservatism may be removed by tightening the bounds or finding alternative derivations, we take a different approach and resort to a computational procedure to establish the same convergence rate. Specifically, our procedure relies on only a few analytic results derived from basic properties of convex functions, whereas the rest is based on computing a numerical certificate that proves the $\bigO(1/K)$ convergence. The computation only involves determining the feasibility of an LMI and can be carried out efficiently.

\section{Convergence Analysis\label{sec:analysis}}

\subsection{Basic properties of convex functions}

The following basic properties of convex functions will be useful later in the paper. A proof of these can be found in standard literature in convex optimization (cf.~\cite{bubeck_convex_2015}). 
\begin{prop}
[Basic properties of convex functions] Suppose $f\colon\mathbb{R}^{n}\to\mathbb{R}$ is convex and differentiable, and $x$ and $y$ are any vectors in $\mathbb{R}^{n}$. Then we have 
\begin{equation}
f(x)\geq f(y)+\nabla f(y)^{T}(x-y).\label{eq:cvx_1st_order}
\end{equation}
Moreover, if $f$ is also $\beta$-smooth, then we have 
\begin{gather}
f(x)\leq f(y)+\nabla f(y)^{T}(x-y)+\frac{\beta}{2}\norm{x-y}^{2},\label{eq:cvx_ub}\\
(\nabla f(x)-\nabla f(y))^{T}(x-y)\geq\frac{1}{\beta}\norm{\nabla f(x)-\nabla f(y)}^{2}.\label{eq:co-coerc}
\end{gather}
\end{prop}
The last property is commonly known as \emph{co-coercivity} of the gradient. 

\subsection{Optimization algorithms as a feedback interconnection\label{subsec:Stability-analysis-of}}

\begin{figure}
\begin{centering}
\includegraphics[width=0.45\columnwidth]{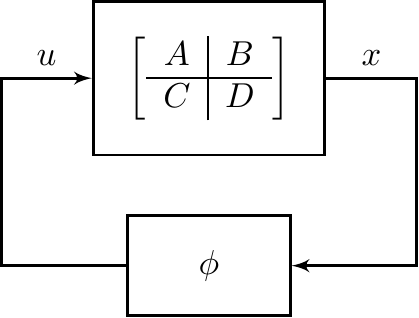}
\par\end{centering}
\caption{Optimization algorithms as a feedback interconnection.\label{fig:fbk_inter}}
\end{figure}

It has been shown in~\cite{lessard_analysis_2016} that many first-order optimization methods can be described as a feedback interconnection (Fig.~\ref{fig:fbk_inter}) of a linear time-invariant (LTI) system and a time-invariant, memoryless, and possibly nonlinear map $\phi$:\begin{subequations}\label{eq:fbk_inter}
\begin{align}
\xi^{k+1} & =A\xi^{k}+Bu^{k},\qquad x^{k}=C\xi^{k}+Du^{k},\label{eq:lti}\\
u^{k} & =\phi(x^{k}),\label{eq:nl_map}
\end{align}
\end{subequations}where $u$ and $x$ are the input and output of the LTI system, respectively. For the distributed gradient tracking algorithm in~(\ref{eq:dgd}), we can choose the state variable $\xi=(\xi^{1},\xi^{2})\deq(x,s-u)$ and the nonlinear map $\phi\colon\mathbb{R}^{n}\to\mathbb{R}^{n}$ as $\phi\coloneqq g$, where $g$ is defined in~(\ref{eq:gk}). Then, it can be shown that the corresponding LTI system in~(\ref{eq:fbk_inter}) is given by
\begin{equation}
\left[\begin{array}{c|c}
A & B\\
\hline C & D
\end{array}\right]=\left[\begin{array}{cc|c}
W & \eta I & -\eta I\\
0 & W & W-I\\
\hline I & 0 & 0
\end{array}\right].\label{eq:lti_dgd}
\end{equation}
The feedback interconnection is well-posed because $D=0$. The dynamics of $\xi^{2}$ imply $\vOne^{T}\xi^{2,k+1}=\vOne^{T}\xi^{2,k}$ for all $k\geq0$. Also, from the initial condition $s^{0}=g(x^{0})=u^{0}$, we have $\xi^{2,0}=0$. Therefore, we have 
\begin{equation}
\vOne^{T}\xi^{2,k}=0,\quad\forall k\geq0.\label{eq:xi2_constr}
\end{equation}
In steady state, we have $\xopt=\vOne\xmin$, $s\opt=0$, and $\uopt=g(\vOne\xmin)$, and therefore we have $\xiopt=(\vOne\xmin,-g(\vOne\xmin))$. For convergence analysis of optimization algorithms, we seek to characterize the behavior of an algorithm under a class of objective functions, and therefore we do not assume to know the exact form of $\phi$. Instead, we characterize $\phi$ through constraints that describe the behavior of its input $x$ and output $u=\phi(x)$. One particularly useful form of constraints is IQCs. The simplest form of IQCs is a \emph{pointwise} IQC, which is a quadratic inequality of the form
\[
\left[\begin{array}{c}
\hat{x}^{k}\\
\hat{u}^{k}
\end{array}\right]^{T}M\left[\begin{array}{c}
\hat{x}^{k}\\
\hat{u}^{k}
\end{array}\right]\geq0
\]
or equivalently 
\begin{align}
 & \psi(\xi^{k},u^{k};M)\deq\nonumber \\
 & \quad\left[\begin{array}{c}
\hat{\xi}^{k}\\
\hat{u}^{k}
\end{array}\right]^{T}\left[\begin{array}{cc}
C & D\\
0 & I
\end{array}\right]^{T}M\left[\begin{array}{cc}
C & D\\
0 & I
\end{array}\right]\left[\begin{array}{c}
\hat{\xi}^{k}\\
\hat{u}^{k}
\end{array}\right]\geq0\label{eq:pwiqc}
\end{align}
for all $k\geq0$. In the above inequalities, the sequences $\hat{\xi}$, $\hat{x}$, and $\hat{u}$ are defined as 
\[
\hat{\xi}^{k}\coloneqq\xi^{k}-\xiopt,\qquad\hat{u}^{k}\coloneqq u^{k}-\uopt,\qquad\hat{x}^{k}\coloneqq x^{k}-\xopt,
\]
and $M$ is a symmetric but possibly indefinite matrix. For $\phi=g$ , it can be shown (Theorem~\ref{thm:main}) that $(\hat{\xi}^{k},\hat{u}^{k})$ is characterized by a pointwise IQC of the form~(\ref{eq:pwiqc}) with 
\[
M=\left[\begin{array}{cc}
0 & \beta I_{n}\\
* & -2I_{n}
\end{array}\right].
\]
We will only use pointwise IQCs in this paper, but we refer the readers to~\cite{lessard_analysis_2016} for more examples of IQCs used in the analysis of optimization algorithms. Note that IQCs should not be used to handle equality constraints of the form $F\hat{\xi}^{k}+G\hat{u}^{k}=0$, which often appear in distributed algorithms. For example, in the distributed gradient tracking algorithm, using~(\ref{eq:xi2_constr}) and the fact $\vOne^{T}\xi^{2,\star}=-\vOne^{T}g(\vOne\xopt_{0})=0$, we can obtain an equality constraint with
\begin{equation}
F=\left[\begin{array}{cc}
0 & \vOne^{T}\end{array}\right],\qquad G=0.\label{eq:FG}
\end{equation}

\subsection{Convergence analysis using dissipativity theory\label{subsec:Bounding-convergence-rate}}

The main mathematical tool we use in analyzing optimization algorithms is dissipativity theory. A (discrete-time) dynamical system with state $\xi$ and input $u$ is called \emph{dissipative} with respect to a \emph{supply rate} $\sr$ if there exists a positive semidefinite function $V$ such that 
\begin{equation}
V(\xi^{k+1})-V(\xi^{k})\leq\sr(\xi^{k},u^{k})\label{eq:dissip_ineq}
\end{equation}
holds for all $k\geq0$ along the system trajectory. The function $V$ is called a \emph{storage function}. It has been shown in~\cite[Th.~7]{hu_dissipativity_2017} that dissipativity theory can be used to establish the $\bigO(1/K)$ convergence of centralized optimization algorithms. We will show that the same result can be modified to prove the $O(1/K)$ convergence of the distributed gradient tracking algorithm.
\begin{prop}
\label{prop:dissip_bound} Suppose there exist a storage function~$V$ and a supply rate $\sr$ that satisfy the inequality~(\ref{eq:dissip_ineq}). In addition, suppose $\sr$ satisfies 
\begin{equation}
\sr(\xi^{k},u^{k})\leq-\frac{1}{n}\sum_{i=1}^{n}\left[f_{0}(x_{i}^{k})-\fopt\right]\label{eq:sr_ineq}
\end{equation}
 for all $k\geq0$. Then we have
\[
\frac{1}{n}\sum_{i=1}^{n}\left[f_{0}(\hat{x}_{i}^{K})-\fopt\right]\leq\frac{V(\xi^{0})}{K+1}.
\]
\end{prop}
\begin{IEEEproof}
Due to convexity, we have 
\[
f_{0}(\hat{x}_{i}^{K})=f_{0}\left(\frac{1}{K+1}\sum_{k=0}^{K}x_{i}^{k}\right)\leq\frac{1}{K+1}\sum_{k=0}^{K}f_{0}(x_{i}^{k}).
\]
Define $\delta^{k}\deq-\frac{1}{n}\sum_{i=1}^{n}\left[f_{0}(x_{i}^{k})-\fopt\right]$. Consider the inequality~(\ref{eq:dissip_ineq}) for $k=0,1,\dots,K$, and sum up all the $K+1$ inequalities to obtain
\[
V(\xi^{k+1})-V(\xi^{0})\leq\sum_{k=0}^{K}\sr(\xi^{k},u^{k})\leq\sum_{k=0}^{K}\delta^{k}.
\]
Then, we have 
\[
\frac{1}{n}\sum_{i=1}^{n}\left[f_{0}(\hat{x}_{i}^{K})-\fopt\right]\leq-\frac{1}{K+1}\sum_{k=0}^{K}\delta^{k}\leq\frac{V(\xi^{0})}{K+1},
\]
where in the last step we have used the fact $V(\xi^{k+1})\geq0$. 
\end{IEEEproof}
In order to successfully find a supply rate $\sr$ that satisfies both~(\ref{eq:dissip_ineq}) and~(\ref{eq:sr_ineq}), it is important to construct a rich enough class of supply rates to begin with. As suggested in~\cite[Rem.~6]{hu_dissipativity_2017}, IQCs can be used to generate a class of supply rates in a way described as follows. One can begin with a base supply rate $\sr_{0}$ that satisfies~(\ref{eq:sr_ineq}). Suppose $(\xi^{k},u^{k})$ can be characterized by a total of $p$ IQCs whose corresponding matrices are given by $\{M_{i}\}_{i=1}^{p}$. Then, it can be seen from~(\ref{eq:pwiqc}) that the supply rate $\sigma\deq\sigma_{0}-\sum_{i=1}^{p}\lambda_{i}\psi(\cdot,\cdot;M_{i})$ also satisfies~(\ref{eq:sr_ineq}) for any $\lambda_{1},\lambda_{2},\dots,\lambda_{p}\geq0$. For satisfying the other inequality~(\ref{eq:dissip_ineq}), when $\sigma_{0}$ is quadratic, i.e., 
\begin{equation}
\sigma_{0}(\xi^{k},u^{k})=\left[\begin{array}{c}
\hat{\xi}^{k}\\
\hat{u}^{k}
\end{array}\right]^{T}S_{0}\left[\begin{array}{c}
\hat{\xi}^{k}\\
\hat{u}^{k}
\end{array}\right]\label{eq:quadratic_sr}
\end{equation}
for some symmetric matrix $S_{0}$, we may choose to restrict our search of the storage function $V$ also to quadratic functions, i.e., $V(\xi^{k})=(\hat{\xi}^{k})^{T}P\hat{\xi}^{k}$, where $P$ is to be determined. Then, as shown in~\cite{hu_dissipativity_2017} and~\cite[Lem.~2]{sundararajan_robust_2017}, the problem of finding $V$ and $\sigma$ that satisfy~(\ref{eq:dissip_ineq}) becomes finding $P\succeq0$ and $\lambda_{i}\geq0$ such that
\begin{align}
 & R^{T}\left\{ \left[\begin{array}{cc}
A^{T}PA-P & A^{T}PB\\
* & B^{T}PB
\end{array}\right]-S_{0}\right.\nonumber \\
 & \quad\left.+\sum_{i=1}^{p}\lambda_{i}\left[\begin{array}{cc}
C & D\\
0 & I
\end{array}\right]^{T}M_{i}\left[\begin{array}{cc}
C & D\\
0 & I
\end{array}\right]\right\} R\preceq0,\label{eq:sr_lmi}
\end{align}
where $R$ is a matrix whose columns span the nullspace of $\left[\begin{array}{cc}
F & G\end{array}\right]$. 

\subsection{Choice of supply rate}

We need to choose the base supply rate~$\sr_{0}$ and IQCs that characterize smooth but non-strongly convex objective functions. Note that such a choice should only depend on properties of the objective function and not on the optimization algorithm. We first state a technical lemma derived from basic properties of smooth convex functions. 
\begin{lem}
\label{lem:smooth_ineq}For all $v,w\in\mathbb{R}$ and $x\in\mathbb{R}^{n}$, we have
\begin{equation}
f_{0}(v)-f_{0}(w)\leq\bar{g}(x)(v-w)+\frac{\beta}{2n}\norm{x-\mathbf{1}v}^{2}.\label{eq:fv_fw}
\end{equation}
where $\bar{g}(x)\coloneqq\frac{1}{n}\sum_{i=1}^{n}\nabla f_{i}(x_{i})=\frac{1}{n}\vOne^{T}g(x).$
\end{lem}
\begin{IEEEproof}
For any $v\in\mathbb{R}$ and $x\in\mathbb{R}^{n}$, because $f_{i}$ is convex, we have from~(\ref{eq:cvx_1st_order})
\begin{align}
f_{0}(v) & \geq\frac{1}{n}\sum_{i=1}^{n}\left[f_{i}(x_{i})+\nabla f_{i}(x_{i})(v-x_{i})\right].\label{eq:fv_geq}
\end{align}
Also, because $f_{i}$ is $\beta$-smooth, we have from~(\ref{eq:cvx_ub})
\begin{equation}
f_{0}(v)\leq\frac{1}{n}\sum_{i=1}^{n}\left[f_{i}(x_{i})+\nabla f_{i}(x_{i})(v-x_{i})+\frac{\beta}{2}\norm{v-x_{i}}^{2}\right].\label{eq:fv_leq}
\end{equation}
For any $v,w\in\mathbb{R}$ and $x\in\mathbb{R}^{n}$, using both~(\ref{eq:fv_geq}) and~(\ref{eq:fv_leq}) and changing $v\to w$ in~(\ref{eq:fv_geq}), we have 
\begin{align*}
f_{0}(v)-f_{0}(w) & \leq\frac{1}{n}\sum_{i=1}^{n}\left[\nabla f_{i}(x_{i})(v-w)+\frac{\beta}{2}\norm{v-x_{i}}^{2}\right]\\
 & =\bar{g}(x)(v-w)+\frac{\beta}{2n}\norm{x-\vOne v}^{2}.
\end{align*}
\end{IEEEproof}
Based on this lemma, we can derive a quadratic supply rate $\sr_{0}$ that satisfies~(\ref{eq:sr_ineq}). This supply rate will be used later in combination with IQCs for numerically proving the $O(1/K)$ convergence of the distributed gradient tracking algorithm.
\begin{prop}
[Choice of supply rate]\label{prop:dgd_zks}Consider a supply rate $\sr_{0}$ of the quadratic form given in~(\ref{eq:quadratic_sr}) with
\begin{align}
 & S_{0}=-\frac{1}{n}\left[\begin{array}{cc}
\beta C^{T}\matdJ C & C^{T}\matJ/2\\
* & 0
\end{array}\right],\label{eq:storage_func}
\end{align}
where $\matJ$ and $\matdJ$ are given by $\matJ\deq\frac{1}{n}\mathbf{1}\mathbf{1}^{T}$ and $\matdJ\deq I-\matJ$, and the matrix $C$ is given in~(\ref{eq:lti_dgd}). We have $\sr_{0}(\xi^{k},u^{k})\leq-\frac{1}{n}\sum_{i=1}^{n}\left[f_{0}(x_{i}^{k})-\fopt\right]$ for all $k\geq0$.
\end{prop}
\begin{IEEEproof}
Using Lemma~\ref{lem:smooth_ineq} and changing $(v,w,x)\to(x_{i}^{k},\xmin,x^{k})$ in~(\ref{eq:fv_fw}), we have
\begin{equation}
f(x_{i}^{k})-\fopt\leq\bar{g}(x^{k})(x_{i}^{k}-\xmin)+\frac{\beta}{2n}\norm{x^{k}-\mathbf{1}x_{i}^{k}}^{2}.\label{eq:ffopt_two_terms}
\end{equation}
Summing up~(\ref{eq:ffopt_two_terms}) over $i=1,2,\dots,n$, we obtain
\begin{align}
 & \sum_{i=1}^{n}\left[f(x_{i}^{k})-\fopt\right]\nonumber \\
 & \quad\leq\bar{g}(x^{k})\vOne^{T}(x^{k}-\xmin)+\frac{\beta}{2n}\sum_{i=1}^{n}\norm{x^{k}-\mathbf{1}x_{i}^{k}}^{2}\nonumber \\
 & \quad=\bar{g}(x^{k})\vOne^{T}(x^{k}-\xmin)+\beta(\hat{x}^{k})^{T}\matdJ x^{k}.\label{eq:prop_sr_eqn1}
\end{align}
Using the facts $\vOne^{T}\uopt=0$ and $\matdJ\xopt=0$, we can rewrite the right-hand side of~(\ref{eq:prop_sr_eqn1}) as
\[
\frac{1}{n}(\hat{u}^{k})^{T}\vOne\vOne\hat{x}^{k}+\beta(\hat{x}^{k})^{T}\matdJ\hat{x}^{k}=-n\left[\begin{array}{c}
\hat{\xi}^{k}\\
\hat{u}^{k}
\end{array}\right]^{T}S_{0}\left[\begin{array}{c}
\hat{\xi}^{k}\\
\hat{u}^{k}
\end{array}\right],
\]
therefore completing the proof. 
\end{IEEEproof}

With the base supply rate $\sr_{0}$, we are ready to present the main result of this paper, which is a numerical procedure for proving the $\bigO(1/K)$ convergence of the distributed gradient tracking algorithm for $\beta$-smooth objective functions. 
\begin{thm}
\label{thm:main}Consider the LMI (\ref{eq:sr_lmi}) for $p=1$, where $(A,B,C)$ is given by (\ref{eq:lti_dgd}), the matrices $F$ and $G$ are given by (\ref{eq:FG}), the matrix $M_{1}$ is given by
\[
M_{1}=\left[\begin{array}{cc}
0 & \beta I_{n}\\
* & -2I_{n}
\end{array}\right],
\]
and $S_{0}$ is given by (\ref{eq:storage_func}). Suppose there exist $P\succeq0$ and $\lambda_{1}\geq0$ such that the LMI (\ref{eq:sr_lmi}) is feasible. Then, the sequence $\{x^{k}\}$ generated by the distributed gradient tracking algorithm (\ref{eq:dgd}) satisfies
\[
\frac{1}{n}\sum_{i=1}^{n}\left[f_{0}(\hat{x}_{i}^{K})-\fopt\right]\leq\frac{V_{0}}{K+1}
\]
for some constant $V_{0}>0$, where $\hat{x}_{i}^{K}=\frac{1}{K+1}\sum_{k=0}^{K}x_{i}^{k}$. 
\end{thm}
\begin{IEEEproof}
Because $f_{i}$ is $\beta$-smooth, we have from (\ref{eq:co-coerc})
\begin{align*}
 & \left[\begin{array}{c}
x_{i}^{k}-\xmin\\
\nabla f_{i}(x_{i}^{k})-\nabla f_{i}(\xmin)
\end{array}\right]^{T}\\
 & \qquad\left[\begin{array}{cc}
0 & \beta\\
* & -2
\end{array}\right]\left[\begin{array}{c}
x_{i}^{k}-\xmin\\
\nabla f_{i}(x_{i}^{k})-\nabla f_{i}(\xmin)
\end{array}\right]\geq0.
\end{align*}
Collecting all the above inequalities for $i=1,2,\dots,n$, we can obtain a pointwise IQC of the form~(\ref{eq:pwiqc}) with the associated matrix given by $M_{1}$. If the LMI~(\ref{eq:sr_lmi}) is feasible, then from Proposition~\ref{prop:dissip_bound} and our choice of the supply rate $\sigma_{0}$ based on Proposition~\ref{prop:dgd_zks}, we have 
\begin{align*}
 & \frac{1}{n}\sum_{i=1}^{n}\left[f_{0}(\hat{x}_{i}^{K})-\fopt\right]\leq\frac{1}{K+1}(\hat{\xi}^{0})^{T}P\hat{\xi}^{0},
\end{align*}
which completes the proof if we define $V_{0}\coloneqq(\hat{\xi}^{0})^{T}P\hat{\xi}^{0}$.
\end{IEEEproof}
Any feasible solution $(P,\lambda_{1})$ serves as a numerical certificate of the $O(1/K)$ convergence of the distributed gradient tracking algorithm. For any given $\beta$, we can combine the LMI feasibility problem with a bisection search on $\eta$ to obtain the maximum allowed step size that ensures convergence.

\section{Discussions\label{sec:discussions}}

We compare our results with previous work on control-theoretic analysis of optimization algorithms for non-strongly convex objective functions and also the original result in~\cite{qu_harnessing_2018} on the convergence rate of the distributed gradient tracking algorithm. In all numerical experiments, the LMI~(\ref{eq:sr_lmi}) was solved using CVX (v2.1, build 1123) in MATLAB with MOSEK (v8.1.0.75) as the optimization solver. 

\subsection{Comparison with previous work on the centralized gradient method}

By viewing centralized optimization as a special case of distributed optimization for $n=1$, we are able to compare our result with previous work on control-theoretic analysis of optimization algorithms. In particular, we compare our result with the one by Hu and Lessard~\cite{hu_dissipativity_2017} and by Fazlyab et al.~\cite{fazlyab_analysis_2018}, both of which include an explicit treatment of non-strongly convex objective functions but only for centralized algorithms. For the centralized gradient descent method, the corresponding LTI system in the feedback interconnection~(\ref{eq:fbk_inter}) is given by
\[
\left[\begin{array}{c|c}
A & B\\
\hline C & D
\end{array}\right]=\left[\begin{array}{c|c}
1 & -\eta\\
\hline 1 & 0
\end{array}\right].
\]
By setting $n=1$, we can obtain the matrix $S_{0}$ associated with the supply rate in Proposition~\ref{prop:dgd_zks} and the matrix $M_{1}$ associated the with IQC in Theorem~\ref{thm:main} as

\begin{gather*}
S_{0}=\left[\begin{array}{cc}
0 & -\frac{1}{2}\\
* & 0
\end{array}\right],\qquad M_{1}=\left[\begin{array}{cc}
0 & \beta\\
* & -2
\end{array}\right].
\end{gather*}
(The equality constraint specified by~(\ref{eq:FG}) is no longer in effect.) For $P\geq0$ and $\lambda_{1}\geq0$, the LMI~(\ref{eq:sr_lmi}) becomes
\[
\left[\begin{array}{cc}
0 & -\eta P+\lambda_{1}\beta+\frac{1}{2}\\
* & \eta^{2}P-2\lambda_{1}
\end{array}\right]\preceq0,
\]
which is feasible if and only if $-\eta P+\lambda_{1}\beta+\frac{1}{2}=0$ and $\eta^{2}P-2\lambda_{1}\leq0$. It can be derived that the step size $\eta$ must satisfy $0<\eta<\frac{2}{\beta}$. This result is more general than the condition $0<\eta\leq\frac{1}{\beta}$ given in~\cite[Sec. 4.1]{hu_dissipativity_2017}, which also uses dissipativity theory but with a different supply rate. 

Although the same condition on $\eta$ is obtained in~\cite[eq.~(4.16b)]{fazlyab_analysis_2018} for the centralized case, we would like to point out that the IQC framework in~\cite{fazlyab_analysis_2018} does not appear to be directly applicable to the case of distributed optimization algorithms. One reason is that many of the derivations therein rely on the fact $\uopt=0$, which holds for centralized optimization but not in a distributed setting, in which we have $\uopt=g(\vOne\xmin)$. Another reason is that the time-varying Lyapunov function used in~\cite{fazlyab_analysis_2018} can be difficult to find. According to~\cite[Th.~3.1]{fazlyab_analysis_2018}, finding such a Lyapunov function is equivalent to solving the following LMI in $P$, $\{a^{k}\}$, and $\{\sigma^{k}\}$: 
\[
M_{0}^{k}(P)+a^{k}M_{1}^{k}+(a^{k+1}-a^{k})M_{2}^{k}+\sigma^{k}M_{3}^{k}\preceq0
\]
for all $k\geq0$, where $M_{0}^{k}$ is linear in $P$, and $M_{i}^{k}$ ($i=1,2,3)$ is constant. For our problem, the matrix $M_{1}^{k}$ must satisfy
\[
\frac{1}{n}\sum_{i=1}^{n}\left[f_{0}(\hat{x}_{i}^{k+1})-f_{0}(\hat{x}_{i}^{k})\right]\leq\left[\begin{array}{c}
\hat{\xi}^{k}\\
\hat{u}^{k}
\end{array}\right]^{T}M_{1}^{k}\left[\begin{array}{c}
\hat{\xi}^{k}\\
\hat{u}^{k}
\end{array}\right].
\]
In order to prove $\bigO(1/K)$ convergence, the sequence $\{a^{k}\}$ must be nonnegative and increase linearly with $k$. This implies that $M_{1}^{k}$ must eventually become negative semidefinite as $k$ grows, and therefore $\sum_{i=1}^{n}\left[f_{0}(\hat{x}_{i}^{k+1})-f_{0}(\hat{x}_{i}^{k})\right]\leq0$ must always hold for large enough $k$. This condition is stronger than the condition~(\ref{eq:sr_ineq}) required by dissipativity theory, which imposes no condition on the decrement of the objective function after one step of update. 

\subsection{Comparison with previous analytic bound}

The $O(1/K)$ convergence of the distributed gradient tracking method was first established in~\cite{qu_harnessing_2018} using an analytic approach. The convergence proof therein uses explicitly the spectral properties of $W$. In comparison, our method relies on only a few basic properties of convex functions, whereas the spectral properties of $W$ are encoded implicitly into the LMI feasibility problem. 

\begin{figure}
\begin{centering}
\includegraphics[width=0.7\columnwidth]{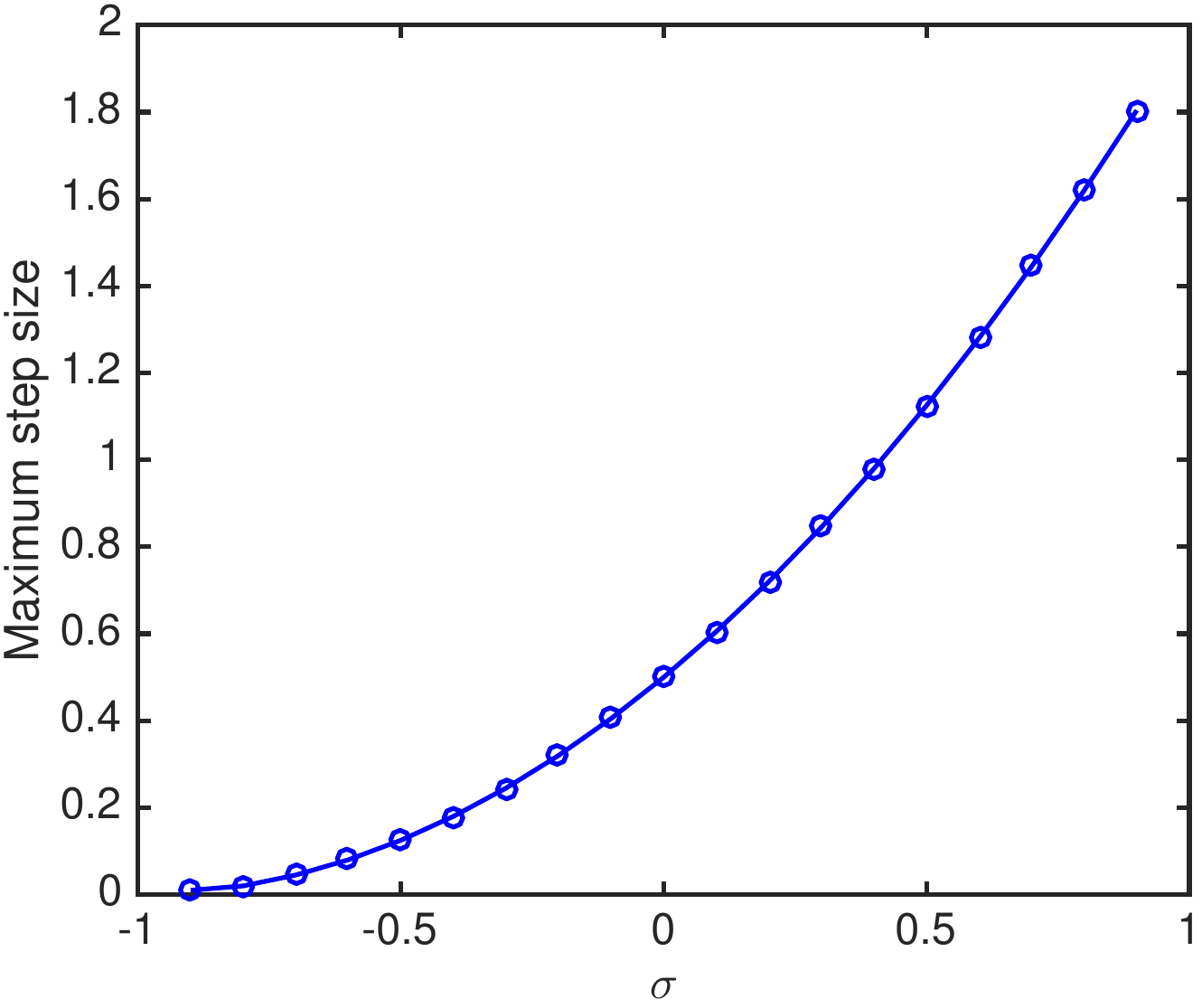}
\par\end{centering}
\caption{Maximum step size as a function of $\sigma$, which defines $W$ in~(\ref{eq:W_def}).\label{fig:max_step_size}}

\end{figure}

Compared to the analytic approach, the computational approach developed in this paper is capable of certifying convergence for a broader range of step sizes $\eta$. We believe that this is because the computational approach does not rely on relaxations for the purpose of maintaining closed-form expressions, which are necessary in the analytic approach. For the purpose of illustration, we conducted a numerical experiment, in which we chose $\beta=1$ and 
\begin{equation}
W=\left[\begin{array}{cc}
\frac{1+\sigma}{2} & \frac{1-\sigma}{2}\\
\frac{1-\sigma}{2} & \frac{1+\sigma}{2}
\end{array}\right],\label{eq:W_def}
\end{equation}
where $\sigma\in(-1,1)$. For different values of $\sigma$, we computed the maximum step size (Fig.~\ref{fig:max_step_size}) that guarantees convergence by solving the LMI~(\ref{eq:sr_ineq}) with a bisection search on $\eta$. In comparison, the result in~\cite{qu_harnessing_2018} only guarantees convergence for $0<\eta\leq\frac{(1-\sigma_{2})^{2}}{160\beta}$, where $\sigma_{2}$ is the second-largest singular value of $W$. Consider, for example, the case when $\sigma=0.5$, in which second-largest singular value of $W$ is given by $\sigma_{2}=|\sigma|=0.5$. The computational approach can ensure convergence with a maximum step size of $1.13$, whereas the analytic approach is more conservative and can only permit a step size up to $\frac{(1-\sigma_{2})^{2}}{160\beta}=1.56\times10^{-3}$. Another important feature of the computational approach is that it is able to take the actual value of $W$ into account when determining the step size, whereas the analytic approach only uses the second-largest singular value $\sigma_{2}$ of $W$. This can make a difference when the value of $W$ is known: for example, both $\sigma=0.5$ and $\sigma=-0.5$ in~(\ref{eq:W_def}) give the same $\sigma_{2}$ but allow different maximum step sizes as shown in Fig.~\ref{fig:max_step_size}.

\subsection{Generalization\label{subsec:Generalizations}}

Although Theorem~\ref{thm:main} is presented for the case of $d=1$, it can be easily generalized to the case of $d>1$ with no modification. When $d>1$, we can replace the definition of $x^{k}$ in~(\ref{eq:xk}) by 
\[
x^{k}\coloneqq\left[\begin{array}{cccc}
(x_{1}^{k})^{T} & (x_{2}^{k})^{T} & \cdots & (x_{n}^{k})^{T}\end{array}\right]^{T}\in\mathbb{R}^{nd}.
\]
Under the new definition of $x^{k}$, it can be shown that the matrices $A$, $B$, $C$, $\{M_{i}\}_{i=1}^{p}$, $S_{0}$, and $R$ in the LMI~(\ref{eq:sr_lmi}) are replaced by their Kronecker products with $I_{d}$, e.g., changing $A\to A\otimes I_{d}$. Then, it follows from properties of the Kronecker product that the new LMI (with Kronecker products) is feasible if and only if the original LMI~(\ref{eq:sr_lmi}) is feasible. Therefore, the same LMI~(\ref{eq:sr_lmi}) can still be used to certify convergence for $d>1$, which is known in existing literature as \emph{lossless dimensionality reduction} (cf.~\cite[Sec.~4.2]{lessard_analysis_2016}).

\section{Conclusions}

We have presented a computational analysis of the $\bigO(1/K)$ convergence of the distributed gradient tracking algorithm in~\cite{qu_harnessing_2018} when the objective function is smooth and convex but not strongly convex. The analysis is built upon existing work on analyzing optimization algorithms using dissipativity theory. We have derived a base supply rate for this setting that involves both non-strongly convex objective functions and distributed optimization. The base supply rate is used to generate a class of supply rates that is rich enough to ensure a successful search of a numerical convergence certificate consisting of a suitable supply rate and storage function. Similar to existing work, the numerical certificate can be automatically generated by solving an LMI feasibility problem. Our method only requires a few analytic derivations from basic properties of convex functions and has been found to prove convergence under a broader range of step sizes than the previous analytic result in~\cite{qu_harnessing_2018}. We expect no difficulty in using the same method for analyzing other first-order distributed optimization algorithms that can be written as the feedback interconnection~(\ref{eq:fbk_inter}) such as EXTRA~\cite{shi_extra_2015} and NIDS~\cite{li_decentralized_2017}. (See~\cite{sundararajan_robust_2017} for state space realizations of these algorithms.)

\bibliographystyle{abbrv}
\bibliography{ACC2019}

\begin{thebibliography}{10}

\bibitem{bach_non-strongly-convex_2013}
F.~Bach and E.~Moulines.
\newblock Non-strongly-convex smooth stochastic approximation with convergence
  rate ${O}(1/n)$.
\newblock In {\em International {Conference} on {Neural} {Information}
  {Processing} {Systems} ({NIPS})}, pages 773--781, 2013.

\bibitem{boyd_distributed_2011}
S.~Boyd, N.~Parikh, E.~Chu, B.~Peleato, and J.~Eckstein.
\newblock Distributed optimization and statistical learning via the alternating
  direction method of multipliers.
\newblock {\em Foundations and Trends in Machine Learning}, 3(1):1--122, 2011.

\bibitem{bubeck_convex_2015}
S.~Bubeck.
\newblock Convex {Optimization}: {Algorithms} and {Complexity}.
\newblock {\em Foundations and Trends in Machine Learning}, 8(3-4):231--357,
  2015.

\bibitem{fazlyab_analysis_2018}
M.~Fazlyab, A.~Ribeiro, M.~Morari, and V.~M. Preciado.
\newblock Analysis of optimization algorithms via integral quadratic
  constraints: {Nonstrongly} convex problems.
\newblock {\em SIAM Journal on Optimization}, 28(3):2654--2689, 2018.

\bibitem{hu_dissipativity_2017}
B.~Hu and L.~Lessard.
\newblock Dissipativity theory for {Nesterov}'s accelerated method.
\newblock In {\em International {Conference} on {Machine} {Learning} ({ICML})},
  pages 1549--1557, 2017.

\bibitem{lessard_analysis_2016}
L.~Lessard, B.~Recht, and A.~Packard.
\newblock Analysis and design of optimization algorithms via integral quadratic
  constraints.
\newblock {\em SIAM Journal on Optimization}, 26(1):57--95, 2016.

\bibitem{li_decentralized_2017}
Z.~Li, W.~Shi, and M.~Yan.
\newblock A decentralized proximal-gradient method with network independent
  step-sizes and separated convergence rates.
\newblock {\em arXiv:1704.07807}, 2017.

\bibitem{megretski_system_1997}
A.~Megretski and A.~Rantzer.
\newblock System analysis via integral quadratic constraints.
\newblock {\em IEEE Transactions on Automatic Control}, 42(6):819--830, 1997.

\bibitem{molzahn_survey_2017}
D.~K. Molzahn, F.~D{\"o}rfler, H.~Sandberg, S.~H. Low, S.~Chakrabarti,
  R.~Baldick, and J.~Lavaei.
\newblock A survey of distributed optimization and control algorithms for
  electric power systems.
\newblock {\em IEEE Transactions on Smart Grid}, 8(6):2941--2962, 2017.

\bibitem{nedic_distributed_2018}
A.~Nedi{\'c} and J.~Liu.
\newblock Distributed optimization for control.
\newblock {\em Annual Review of Control, Robotics, and Autonomous Systems},
  1(1):77--103, 2018.

\bibitem{nedic_achieving_2017}
A.~Nedi{\'c}, A.~Olshevsky, and W.~Shi.
\newblock Achieving geometric convergence for distributed optimization over
  time-varying graphs.
\newblock {\em SIAM Journal on Optimization}, 27(4):2597--2633, 2017.

\bibitem{parrilo_structured_2000}
P.~A. Parrilo.
\newblock {\em Structured semidefinite programs and semialgebraic geometry
  methods in robustness and optimization}.
\newblock {PhD} {Thesis}, California Institute of Technology, 2000.

\bibitem{qu_harnessing_2018}
G.~Qu and N.~Li.
\newblock Harnessing smoothness to accelerate distributed optimization.
\newblock {\em IEEE Transactions on Control of Network Systems},
  5(3):1245--1260, 2018.

\bibitem{rabbat_distributed_2004}
M.~Rabbat and R.~Nowak.
\newblock Distributed optimization in sensor networks.
\newblock In {\em International {Symposium} on {Information} {Processing} in
  {Sensor} {Networks} ({IPSN})}, pages 20--27, 2004.

\bibitem{shi_extra_2015}
W.~Shi, Q.~Ling, G.~Wu, and W.~Yin.
\newblock {EXTRA}: {An} exact first-order algorithm for decentralized consensus
  optimization.
\newblock {\em SIAM Journal on Optimization}, 25(2):944--966, 2015.

\bibitem{sundararajan_robust_2017}
A.~Sundararajan, B.~Hu, and L.~Lessard.
\newblock Robust convergence analysis of distributed optimization algorithms.
\newblock In {\em Annual {Allerton} {Conference} on {Communication}, {Control},
  and {Computing}}, 2017.

\end{thebibliography}

\end{document}